\documentclass[12pt]{article}
\usepackage{bbm}
\usepackage{mathrsfs}
\usepackage{amssymb,amsmath}
\def\D{\Delta}
\def\UU{{\mathcal U}}
\def\XX{{\mathcal X}}
\def\HH{{\mathcal H}}

\def\WW{{\mathcal W}}
\def\LL{\mathcal L}

\def\e{\epsilon}

\def\cl{\centerline}

\def\nz{\mathbbm{n}}
\def\cz{\mathbbm{c}}
\def\bz{\mathbbm{b}}

\def\vs{\vspace*}

\def\DD{{\mathcal D}}

\def\ni{\noindent}

\def\ptl{\hbar}

\def\WW{\mathcal {W}}

\def\Z{\mathbb{Z}}

\def\C{\mathbb{C}}
\def\CC{\mathcal {C}}
\def\QED{\hfill$\Box$}
\numberwithin{equation}{section}
\newtheorem{theo}{Theorem}[section]
\newtheorem{defi}[theo]{Definition}

\newtheorem{lemm}[theo]{Lemma}

\begin{document}

\cl{{\large \bf Quantizations of the $W$ Algebra $W(2,2)$}\footnote
{Supported by NSF grants 10471091, 10671027 of China, ``One Hundred
Talents Program'' from University of Science and Technology of China.\\[2pt]
\indent Corresponding E-mail: sd\_junbo@163.com} } \vs{6pt}

\cl{Junbo Li$^{*,\dag)}$, Yucai Su$^{\ddag)}$} \cl{\small
$^{*)}$Department of Mathematics, Shanghai Jiao Tong University,
 Shanghai, 200240, China}
\cl{\small $^{\dag)}$ Department of Mathematics, Changshu Institute
of Technology, Changshu 215500, China}  \cl{\small
$^{\ddag)}$Department of Mathematics, University of Science and
Technology of China Hefei 230026, China} \cl{\small E-mail:
sd\_junbo@163.com, ycsu@ustc.edu.cn} \vs{6pt}

\noindent{\small{\bf Abstract.} We quantize the $W$-algebra
$W(2,2)$, whose Verma modules, Harish-Chandra modules, irreducible
weight modules and Lie bialgebra structures have been investigated
and determined in a series of papers recently.

  }
\vs{12pt}

\cl{\bf\S1. \
Introduction}\setcounter{section}{1}\setcounter{equation}{0}

It is interesting to construct new quantum groups by quantizing Lie
bialgebras to some authors (see \cite{D2}--\cite{HW},
\cite{LS2,SSW}). Originally, Witt and Virasoro type algebras were
quantized in \cite{G}, whose Lie bialgebra structure were presented
in \cite{T} and further classified in \cite{NT}. Afterwards, the Lie
bialgebras of generalized Witt type, generalized Virasoro-like type
and Block type were considered in \cite{SS}, \cite{WSS} and
\cite{LSX} respectively, which were quantized in \cite{HW},
\cite{SSW} and \cite{LS2}.

In the present paper we shall quantize the $W$ algebra $W(2,2)$,
introduced by Zhang and Dong in \cite{ZD}, whose Lie bialgebra
structures have been proved to be triangular coboundary in
\cite{LS3} by the authors. The algebra $\WW$ considered in this
paper is an infinite-dimensional Lie algebra with a $\C$-basis
$\{\,L_n,\,W_n\,|\,n\in \Z\,\}$ and the following Lie brackets
(\,other components vanishing):
\begin{eqnarray}
&&[L_m,L_n]=(m-n)L_{m+n},\ \ \,[L_m,W_n]=(m-n)W_{m+n}.\label{BLie02}
\end{eqnarray}
The highest weight modules of $\WW$ were investigated in \cite{ZD},
which produce a new class of irrational vertex operator algebras.
Later on, the irreducible weight modules and indecomposable modules
on $\WW$ were considered in \cite{LGZ} and \cite{LZ}.

Let $\mathcal{A}$ be a unitary $\C$-algebra. For any $y \in
\mathcal{A},\,\bz\in\C,\,i\in \Z$, set
\begin{eqnarray*}
&&y_\bz^{<i>}=\mbox{$\prod\limits_{k=0}^{i-1}$}(y+\bz+k),\ \ \
y_\bz^{[i]}=\mbox{$\prod\limits_{k=0}^{i-1}$}(y+\bz-k),\ \ \
\Big(\!\!\begin{array}{c}\bz\\ i
\end{array}\!\!\Big)=\frac{1}{i\,!}\mbox{$\prod\limits_{k=0}^{i-1}$}(\bz-k).
\end{eqnarray*}

\begin{lemm}\label{lemma1.1}{\rm (\cite{G})}
For any $y\in\mathcal{A}$, $\bz,\cz\in \C$, $m,n,k\in \Z$, one has
\begin{eqnarray}
&&\!\!\!\!\!\!\!\!y_\bz^{<m + n>}=y_\bz^{<m>}y_{\bz+m}^{<n>},\ \ \
y_\bz^{[m+n]}=y_\bz^{[m]}y_{\bz-m}^{[n]},
\ \ \ y_\bz^{[m]}=y_{\bz-m+1}^{<m>},\label{eqnayam+n11}\\
&&\!\!\!\!\!\!\!\!\mbox{$\sum\limits_{m+n =
k}$}\frac{(-1)^n}{m\,!\,n\,!}y_\bz^{[m]}y_\cz^{<n>}=\Big(\!\!\begin{array}{c}\bz-\cz
\\ k\end{array}\!\!\Big),\ \
\mbox{$\sum\limits_{m+n=k}$}\frac{(-1)^n
}{m\,!\,n\,!}y_\bz^{[m]}y_{\cz-m}^{[n]}=\Big(\!\!\begin{array}{c}\bz-\cz+k-1
\\ k\end{array}\!\!\Big).\label{eqnayam62+n11}
\end{eqnarray}
\end{lemm}
\begin{defi}\label{Defin1.2}{\rm (\cite{D1,D2})}
Let ($\mathcal {H},\mu,\tau,\D,\e,\mathcal {S}$) be a quantized
enveloping algebra satisfying $\mathcal {H}/t\mathcal {H}\cong
\UU(\LL)$, where $\LL$ is a Lie algebra and $t$ a deformation
parameter. An invertible element $\DD\in\HH\otimes\HH$ is called a
Drinfeld's twisting if it satisfies the following conditions
\begin{eqnarray}
&&(\DD\otimes 1)(\D\otimes Id)(\DD)=(1\otimes \DD)(Id
\otimes \D)(\DD),\label{eq1.4abc}\\
&&(\e\otimes Id )(\DD)=1\otimes1=(Id\otimes
\e)(\DD).\label{eq1.5abc}
\end{eqnarray}
\end{defi}
The following well-known theorem is due to Drinfeld (\,see
\cite{D1,D2}).
\begin{theo}\label{theo1.3}
Let ($\HH,\mu,\tau,\D^o,\e,\mathcal {S}^o$) be a Hopf algebra over a
commutative ring and $\DD$ a Drinfeld's element of $\HH\otimes\HH$.
Then\\
(1) $\mathcal{U}=\mu\cdot(\mathcal {S}^o\otimes Id)(\DD)$ is
invertible
in $\HH$ with $\mathcal{U}^{-1}=\mu(Id\otimes\mathcal {S}^o)(\DD^{-1})$.\\
(2) The algebra $(\HH,\mu,\tau,\D,\e,\mathcal {S})$ is a new Hopf
algebra if we define
\begin{eqnarray*}
\D (h)=\DD\D^o(h)\DD^{-1},\ \ \mathcal {S}(h)
=\mathcal{U}^{-1}\mathcal{S}^o(h)\mathcal{U}.
\end{eqnarray*}
\end{theo}

Denote by $\UU(\WW)$ the universal enveloping algebra of $\mathcal
{W}$ and $(\UU(\WW),\mu,\tau,\D^o,\e,\mathcal {S}^o)$ its natural
Hopf algebra structure. Then for any $x\in\mathcal {W}$,
\begin{eqnarray*}
&&\D^o(x)=x\otimes1+1\otimes x,\,\,\mathcal {S}^o(x)=-x,\,\,\e(x)=0.
\end{eqnarray*}
For any $n_0\in\Z^*$, denote by $\XX=L_{n_0}$ or $\XX=W_{n_0}$ and
$\hbar=-\frac{L_0}{n_0}$, then $[\,\hbar,\XX\,]=\XX$. For any
$\bz\in\C$, denote
\begin{eqnarray}
&&\CC_\bz=\mbox{$\sum\limits_{k=0}^{\infty}$}
\frac{1}{k\,!}\hbar_\bz^{<k>}\otimes\XX^{k}t^k, \ \
\mathcal{U}_\bz=\mu\cdot(\mathcal {S}^o\otimes Id)(\CC_\bz),\label{eqa1.71}\\
&&\DD_\bz=\mbox{$\sum\limits_{k =0}^{\infty }$}
\frac{(-1)^{k}}{k\,!}\hbar_\bz^{[k]}\otimes\XX^k t^k,\ \
\mathcal{V}_\bz=\mu\cdot(Id\otimes\mathcal
{S}^o)(\DD_\bz).\label{eqa1.81}
\end{eqnarray}
Since $\mathcal {S}^o(\hbar_\bz^{<k> })=(-1 )^k \hbar_{-\bz}^{[k] }$
and $\mathcal {S}^o\XX^{k}=(-1)^{k}\XX^k$, one has
\begin{eqnarray}
&&\mathcal{U}_\bz= \mbox{$\sum\limits_{k=0}^{\infty}$}
\frac{(-1)^k}{k\,!} \hbar_{-\bz}^{[k]}\XX^k t^k ,\ \ \
\mathcal{V}_\bz=\mbox{$\sum\limits_{k=0} ^{\infty}$}
\frac{1}{k\,!}\hbar_{\bz}^{[k]}\XX^k t^k \label{va00}.
\end{eqnarray}
For convenience, we shall use the following abbreviations:
\begin{eqnarray}
\hbar_0 ^{<i>}:=\hbar^{<i>},\ \ \hbar_0^{[i]}:=\hbar^{[i]},\ \
\CC_0:=\CC,\ \ \mathcal{D}_0:=\mathcal{D},\ \
\mathcal{U}_0:=\mathcal{U}, \ \ \mathcal{V}_0:=\mathcal{V}.
\end{eqnarray}
The main results of this paper can be formulated as the following
two theorems.
\begin{theo}\label{theo1.4}
If $\XX=L_{n_0}$, then there exists a Hopf structure $(\UU (\mathcal
{W})[[t]],\mu,\tau,\D,\mathcal {S},\e)$ on $\UU (\mathcal {W})[[t]]$
over $\C[[t]]$, such that $\UU(\mathcal {W})[[t]]/t\UU(\mathcal
{W})[[t]]=\UU (\mathcal {W})$, which preserves the product and the
counit of $\UU (\mathcal {W})[[t]]$, while the coproduct and
antipode are determined by
\begin{eqnarray}
\mathcal {S}(L_{n})\!\!\!&=&\!\!\! -(1-\XX
t)^{-\mathbbm{n}}\mbox{$\sum\limits_{k=0}^{\infty}$}b_kL_{n+kn_0}
\hbar_{k}^{<k>}t^k,\\
\mathcal {S}(W_{n})\!\!\!&=&\!\!\! -(1-\XX
t)^{-\nz}\mbox{$\sum\limits_{k=0}^{\infty}$}b_kW_{n+kn_0}
\hbar_{k}^{<k>}t^k,\\
\D(L_n)\!\!\!&=&\!\!\!L_n\otimes (1-\XX t
)^{\nz}+\mbox{$\sum\limits_{k=0}^{\infty}(-1
)^{k}$}b_k\ptl^{<k>}\otimes(1-\XX t)^{-k}L_{n+kn_0}t^k,\\
\D(W_n)\!\!\!&=&\!\!\!W_n\otimes (1-\XX t
)^{\nz}+\mbox{$\sum\limits_{k=0}^{\infty}(-1
)^{k}$}b_k\ptl^{<k>}\otimes(1-\XX t)^{-k}W_{n+kn_0}t^k,
\end{eqnarray}
where $\mathbbm{n}=-\frac{n}{n_0},\,b_k =
\frac{1}{k\,!}\mbox{$\prod\limits_{n=0}^{k-1}$}\big((1-p)n_0-n\big),\
b_0=1$.
\end{theo}
Conveniently, we use the same notations in the following theorem as
those in Theorem \ref{theo1.4}.
\begin{theo}\label{theo1.5}
If $\XX=W_{n_0}$, then there exists another Hopf structure $(\UU
(\mathcal {W})[[t]],\mu,\tau,\D,S,\e)$ on $\UU (\mathcal {W})[[t]]$
over $\C[[t]]$, such that $\UU(\mathcal {W})[[t]]/t\UU(\mathcal
{W})[[t]]=\UU (\mathcal {W})$, which preserves the product and the
counit of $\UU (\mathcal {W})[[t]]$, while the coproduct and
antipode are determined by
\begin{eqnarray}
\mathcal {S}(L_{n})\!\!\!&=&\!\!\! -(1-\XX
t)^{-\mathbbm{n}}\mbox{$\sum\limits_{k=0}^{\infty}$}b_kW_{n+kn_0}
\hbar_{k}^{<k>}t^k,\\
\D(L_n)\!\!\!\!&=&\!\!\!L_n\otimes (1-\XX t
)^{\nz}\!+\!\mbox{$\sum\limits_{k=0}^{\infty}(-1
)^{k}$}b_k\ptl^{<k>}\otimes(1-\XX t)^{-k}W_{n+kn_0}t^k\,,\\
\mathcal {S}(W_{n})\!\!\!&=&\!\!\! -(1-\XX t)^{-\nz}W_{n},\ \ \
\D(W_n)=W_n\otimes (1-\XX t )^{\nz}+1\otimes W_{n}\,.
\end{eqnarray}
\end{theo}

Throughout the paper we denote by $\Z_+$ the set of all nonnegative
integers and $\C$ the set of all complex numbers.

\vskip16pt

\cl{\bf\S2. \ Proofs of the main
results}\setcounter{section}{2}\setcounter{theo}{0}\setcounter{equation}{0}

\vskip6pt

The proofs of Theorem \ref{theo1.4} and Theorem \ref{theo1.5} are
based on a series of technical lemmas, some of which were originally
developed in \cite{G,HW} and ever employed in \cite{LS2,SSW}. For
completeness, we shall still prove them in detail under our
background.
\begin{lemm}\label{lemm2.1}
If $\XX=L_{n_0}$ or $\XX=W_{n_0}$, then for any
$L_{n},L_{m},W_n,W_m\in\mathcal {W}$, one always has
\begin{eqnarray}
&&\XX^k\hbar_\bz^{[i]}=\hbar_{\bz-k }^{[i]}\XX^{k},\
\,\XX^k\hbar_\bz^{<i>}=\hbar_{\bz-k
}^{<i>}\XX^{k},\ \,L_n\hbar_\bz^{[i]}=\hbar_{\bz-\nz}^{[i]}L_n,\label{pxckm}\\
&&L_n\hbar_\bz^{<i>}=\hbar_{\bz-\nz}^{<i>}L_n,\ \,
W_n\hbar_\bz^{[i]}=\hbar_{\bz-\nz}^{[i]}W_n,\ \,
W_n\hbar_\bz^{<i>}=\hbar_{\bz-\nz}^{<i>}W_n,\label{xbjpamj}\\
&&L_n(L_m)^i=\mbox{$\sum\limits_{k=0}^{i}(-1)^{k}$}\Big(\!\!\begin{array}{c}i\\
k \end{array}\!\!\Big)
\mbox{$\prod\limits_{p=0}^{k-1}$}\big((1-p)m-n\big)
(L_m)^{i-k}L_{n+km},\label{xbjxckm}\\
&&L_n(W_m)^i=\mbox{$\sum\limits_{k=0}^{i}(-1)^{k}$}\Big(\!\!\begin{array}{c}i\\
k \end{array}\!\!\Big)
\mbox{$\prod\limits_{p=0}^{k-1}$}\big((1-p)m-n\big)
(W_m)^{i-k}W_{n+km},\label{0801291}\\
&&W_n(L_m)^i=\mbox{$\sum\limits_{k=0}^{i}(-1)^{k}$}\Big(\!\!\begin{array}{c}i\\
k \end{array}\!\!\Big)
\mbox{$\prod\limits_{p=0}^{k-1}$}\big((1-p)m-n\big)
(L_m)^{i-k}W_{n+km},\label{0801292}
\end{eqnarray}
where $m,n\in\Z$, $\nz=-\frac{n}{n_0}$, $i,k\in\Z_+$ and\,
$\bz\in\C$.
\end{lemm}
\ni{\it Proof}\ \ \ From $[\hbar,L_n]=\nz L_n$ and
$[\,\hbar,W_n]=\nz W_n$, one has $L_n\hbar=(\hbar-\nz)L_n$ and
$W_n\hbar=(\hbar-\nz)W_n$, which imply the case $i=1$ of
(\ref{xbjpamj}). We will use induction on $i$ to prove
(\ref{xbjpamj}). Suppose all equations of (\ref{xbjpamj}) hold for
$i$. Then for the case $i+1$, one has
\begin{eqnarray*}
&&L_n\hbar_\bz^{[i+1]}=L_n\hbar_\bz^{[i]}(\hbar+\bz-i)
=\hbar_{\bz-\nz}^{[i]}L_n(\hbar+\bz-i)=\hbar_{\bz-\nz}^{[i+1]}L_n,\\
&&L_n\hbar_\bz^{<i+1>}=L_n\hbar_\bz^{<i>}(\hbar+\bz+i)
=\hbar_{\bz-\nz}^{<i>}L_n(\hbar+\bz+i)=\hbar_{\bz-\nz}^{<i+1>}L_n,\\
&&W_n\hbar_\bz^{[i+1]}=W_n\hbar_\bz^{[i]}(\hbar+\bz-i)
=\hbar_{\bz-\nz}^{[i]}W_n(\hbar+\bz-i)=\hbar_{\bz-\nz}^{[i+1]}W_n,\\
&&L_n\hbar_\bz^{<i+1>}=W_n\hbar_\bz^{<i>}(\hbar+\bz+i)
=\hbar_{\bz-\nz}^{<i>}W_n(\hbar+\bz+i)=\hbar_{\bz-\nz}^{<i+1>}W_n.
\end{eqnarray*}
Then (\ref{xbjpamj}) follows. Similarly, one can obtain
(\ref{pxckm}) from using the induction on $i$ and $k$. Noting that
\begin{eqnarray*}
&&\!\!\!\!\!\!\!\!L_n(L_m)^i\!=\!\mbox{$\sum\limits_{k=0}^{i}(-1)^{k}$}\Big(\!\!
\begin{array}{c}i\\k\end{array}\!\!\Big)(L_m)^{i-k}\big({\rm
ad}L_m\big)^k(L_n),\,\big({\rm ad}L_m\big)^k(L_n)
\!=\!\mbox{$\prod\limits_{p=0}^{k-1}$}\big((1-p)m-n\big)L_{n+km},\\
&&\!\!\!\!\!\!\!\!L_n(W_m)^i\!=\!\mbox{$\sum\limits_{k=0}^{i}(-1)^{k}$}\Big(\!\!
\begin{array}{c}i\\k\end{array}\!\!\Big)(W_m)^{i-k}\big({\rm
ad}W_m\big)^k(L_n),\,\big({\rm ad}W_m\big)^k(L_n)
\!=\!\mbox{$\prod\limits_{p=0}^{k-1}$}\big((1-p)m-n\big)W_{n+km},\\
&&\!\!\!\!\!\!\!\!W_n(L_m)^i\!=\!\mbox{$\sum\limits_{k=0}^{i}(-1)^{k}$}\Big(\!\!
\begin{array}{c}i\\k\end{array}\!\!\Big)(L_m)^{i-k}\big({\rm
ad}L_m\big)^k(W_n),\, \big({\rm ad}L_m\big)^k(W_n)
\!=\!\mbox{$\prod\limits_{p=0}^{k-1}$}\big((1-p)m-n\big)W_{n+km},
\end{eqnarray*}
Then (\ref{xbjxckm})--(\ref{0801292}) follows.\QED

\begin{lemm}\label{lemm2.2} Whether $\XX=L_{n_0}$ or $\XX=W_{n_0}$, for any $\bz,\cz\in\C$, we have
\begin{eqnarray}
\DD_\bz\CC_\cz=1\otimes(1-\XX t)^{\bz-\cz},\ \,\mathcal{V}_\bz
\mathcal{U}_\cz=(1-\XX t)^{-\bz-\cz}.
\end{eqnarray}
Therefore, the elements
$\DD_\bz,\CC_\bz,\mathcal{U}_\bz,\mathcal{V}_\bz$ are invertible
with
$\DD_\bz^{-1}=\CC_{\bz},\,\mathcal{U}_\bz^{-1}=\mathcal{V}_{-\bz}$.
\end{lemm}
\ni{\it Proof}\ \ \ Using (\ref{eqnayam62+n11}), (\ref{eqa1.71}) and
(\ref{eqa1.81}), one has
\begin{eqnarray*}
\DD_\bz\CC_\cz\!\!\!&=&\!\!\!\big(\mbox{$\sum\limits_{i=0}^{\infty
}$} \frac{(-1)^{i}}{i\,!}\hbar_\bz^{[i]}\otimes\XX^i t^i\big)\cdot
\big(\mbox{$\sum\limits_{j=0}^{\infty}$}
\frac{1}{j\,!}\hbar_\cz^{<j>}\otimes\XX^{j}t^j\big)\\
\!\!\!&=&\!\!\!\mbox{$\sum\limits_{i,j=0}^{\infty }$}
\frac{(-1)^{i}}{i\,!j\,!}\hbar_\bz^{[i]}\hbar_\cz^{<j>}\otimes\XX^{i+j}
t^{i+j}=\mbox{$\sum\limits_{k=0}^{\infty }$}(-1)^{k}
\Big(\!\!\begin{array}{c} \bz-\cz\\ k
\end{array}\!\!\Big)\otimes\XX^{k}t^{k}=1\otimes(1-\XX
t)^{\bz-\cz}.
\end{eqnarray*}
Using (\ref{eqnayam62+n11}), (\ref{va00}) and (\ref{pxckm}), one has

\begin{eqnarray*}
\mathcal{V}_\bz
\mathcal{U}_\cz\!\!\!&=&\!\!\!\big(\mbox{$\sum\limits_{i=0}
^{\infty}$} \frac{1}{i\,!}\hbar_{\bz}^{[i]}\XX^i t^i\big)\cdot
\big(\mbox{$\sum\limits_{j=0}^{\infty}$} \frac{(-1)^j}{j\,!}
\hbar_{-\cz}^{[j]}\XX^j t^j\big)=\mbox{$\sum\limits_{i,j=0}^{\infty
}$}
\frac{(-1)^{j}}{i\,!j\,!}\hbar_\bz^{[i]}\hbar_{-\cz-i}^{[j]}\XX^{i+j}t^{i+j}\\
\!\!\!&=&\!\!\!\mbox{$\sum\limits_{k=0}^{\infty }$}\,
\mbox{$\sum\limits_{i+j=k}$}
\frac{(-1)^{j}}{i\,!j\,!}\hbar_\bz^{[i]}\hbar_{-\cz-i}^{[j]}\XX^{k}t^{k}
=\mbox{$\sum\limits_{k=0}^{\infty
}$}\, \Big(\!\!\begin{array}{c}\bz+\cz+k-1
\\ k\end{array}\!\!\Big)\XX^{k}t^{k}
=(1-\XX t)^{-\bz-\cz}.
\end{eqnarray*}
Then this lemma follows.\QED

\begin{lemm}\label{lemm2.3}
For any $\bz \in \C,\,i\in\Z_+$, one can write
\begin{eqnarray}\label{eqDppt1}
&&\D^o\hbar^{[i]}=\mbox{$\sum\limits_{k=0}^{i}$}
\Big(\!\!\begin{array}{c}i\\ k
\end{array}\!\!\Big)\hbar_{-\bz}^{[k]} \otimes\hbar_\bz^{[i-k]} .
\end{eqnarray}
In particular, one has\,\,
$\D^o\hbar^{[i]}=\mbox{$\sum\limits_{k=0}^{i}$}
\Big(\!\!\begin{array}{c}i\\ k
\end{array}\!\!\Big)\hbar^{[k]}\otimes\hbar^{[i-k]}$.
\end{lemm}

\ni{\it Proof}\ \ \ We will use induction on $i$. The case of $i=1$
follows from the formula $\D^o(\hbar)=\hbar\otimes1+1\otimes\hbar$.
Suppose (\ref{eqDppt1}) holds for $i$. As for the case $i+1$, one
has
\begin{eqnarray*}
&&\!\!\!\!\!\!\D^o(\hbar^{[i+1]})=\D^o(\hbar^{[i]})\D^o(\hbar-i)\\
&\!\!\!=\!\!\!&\big(\mbox{$\sum\limits_{k=0}^i$}\Big(\!\!\begin{array}{c}i\\
k\end{array}\!\!\Big)\hbar_{-\bz}^{[k]}\otimes\hbar_\bz^{[i-k]}\big)
\big((\hbar-\bz-i)\otimes1+1\otimes(\hbar+\bz-i)+i(1\otimes1)\big)\\
&\!\!\!=\!\!\!&\big(\mbox{$\sum\limits_{k=1}^{i-1}$}
\Big(\!\!\begin{array}{c}i\\ k \end{array}\!\!\Big)
\hbar_{-\bz}^{[k]}\otimes\hbar_{\bz}^{[i-k]}\big)
\big((\hbar-\bz-i)\otimes1+1\otimes(\hbar+\bz-i)\big)
+\hbar_{-\bz}^{[i]}\otimes(\hbar+\bz-i)\\
&&+i\Big(\mbox{$\sum\limits_{k=0}^{i}$}\Big(\!\!\begin{array}{c}i\\
k\end{array}\!\!\Big)\hbar_{-\bz}^{[k]}\otimes\hbar_{\bz}^{[i-k]}\Big)
+\big(1\otimes\hbar_\bz^{[i+1]}+\hbar_{-\bz}^{[i+1]}\otimes1\big)
+(\hbar-\bz-i)\otimes\hbar_\bz^{[i]}\\
&\!\!\!=\!\!\!&
1\otimes\hbar_\bz^{[i+1]}+\hbar_{-\bz}^{[i+1]}\otimes1+i
\Big(\mbox{$\sum\limits_{k=1}^{i-1}$}\Big(\!\!\begin{array}{c}i\\ k
\end{array}\!\!\Big)\hbar_{-\bz}^{[k]}\otimes\hbar_{\bz}^{[i-k]}\Big)
+(\hbar-\bz)\otimes\hbar_\bz^{[i]}\\
&&+\hbar_{-\bz}^{[i]}\otimes(\hbar+\bz)+
\mbox{$\sum\limits_{k=1}^{i-1}$}\Big(\!\!\begin{array}{c}i\\ k
\end{array}\!\!\Big)\hbar_{-\bz}^{[k+1]}\otimes
\hbar_\bz^{[i-k]}+\mbox{$\sum\limits_{k=1}^{i-1}$}(k-i)
\Big(\!\!\begin{array}{c}i\\ k \end{array}\!\!\Big)
\hbar_{-\bz}^{[k]}\otimes\hbar_\bz^{[i-k]}\\
&&+\mbox{$\sum\limits_{k=1}^{i-1}$}\Big(\!\!\begin{array}{c}i\\ k
\end{array}\!\!\Big)\hbar_{-\bz}^{[k]}\otimes
\hbar_\bz^{[i-k+1]}+\mbox{$\sum\limits_{k=1}^{i-1}(-k)$}
\Big(\!\!\begin{array}{c}i\\ k \end{array}\!\!\Big)
\hbar_{-\bz}^{[k]}\otimes\hbar_\bz^{[i-k]}\\
&\!\!\!=\!\!\!&\big(1\otimes\hbar_\bz^{[i+1]}+\hbar_{-\bz}^{[i+1]}\otimes1\big)
+\Big(\mbox{$\sum\limits_{k=1}^{i-1}$}\Big(\!\!\begin{array}{c}i\\ k
\end{array}\!\!\Big)\hbar_{-\bz}^{[k]}\otimes
\hbar_\bz^{[i-k+1]}+\hbar_{-\bz}^{[i]}\otimes(\hbar+\bz)\Big)\\
&&+\Big((\hbar-\bz)\otimes\hbar_\bz^{[i]}+
\mbox{$\sum\limits_{k=1}^{i-1}$}\Big(\!\!\begin{array}{c}i\\ k
\end{array}\!\!\Big)
\hbar_{-\bz}^{[k+1]}\otimes \hbar_\bz^{[i-k]}\Big)\\
&\!\!\!=\!\!\!&1\otimes\hbar_\bz^{[i+1]}+\hbar_{-\bz}^{[i+1]}\otimes1+
\mbox{$\sum\limits_{k=1}^{i}$}\Big(\Big(\!\!\begin{array}{c}i\\ k-1
\end{array}\!\!\Big)+\Big(\!\!\begin{array}{c}i\\ k \end{array}\!\!\Big)\Big)
\hbar_{-\bz}^{[k]}\otimes\hbar_\bz^{[i+1-k]}\\
&\!\!\!=\!\!\!&\mbox{$\sum\limits_{k=0}^{i+1}$}\Big(\!\!\begin{array}{c}i+1\\
k\end{array}\!\!\Big)\hbar_{-\bz}^{[k]}\otimes\hbar_\bz^{[i+1-k]}.
\end{eqnarray*}
Then this lemma follows.\QED

\begin{lemm}\label{lemma2.4}
$\DD=\sum_{k=0}^{\infty}\frac{(-1)^k}{k\,!}\hbar^{[k]} \otimes
\XX^kt^k$ is a Drinfeld's twist element of $\UU (\mathcal
{W})[[t]]$, i.e., $\DD$ satisfies (\ref{eq1.4abc}) and
(\ref{eq1.5abc}), no matter when $\XX=L_{n_0}$ or $\XX=W_{n_0}$.
\end{lemm}
\ni{\it Proof}\ \ \ Firstly, we check that $\DD$ satisfies
(\ref{eq1.4abc}). Using Lemma \ref{lemm2.3}, (\ref{eqnayam+n11}) and
(\ref{pxckm}), one has
\begin{eqnarray*}
(\DD\otimes 1)(\D^o\otimes Id)(\DD)
&\!\!\!=\!\!\!&(\mbox{$\sum\limits_{i=0}^{\infty}$}\frac{(-1)^i}{i!}\hbar^{[i]}
\otimes\XX^it^i\otimes1)(\D^o\otimes
Id)(\mbox{$\sum\limits_{j=0}^{\infty}$}
\frac{(-1)^j }{j!}\hbar^{[j]}\otimes\XX^jt^j)\\
&\!\!\!=\!\!\!&(\mbox{$\sum\limits_{i=0}^{\infty}$}\frac{(-1)^i}{i!}\hbar^{[i]}
\otimes
\XX^it^i\otimes1)\!\cdot\!(\mbox{$\sum\limits_{j=0}^{\infty}$}\frac{(-1)^j}{j!}
\mbox{$\sum\limits_{k=0}^j$}\Big(\!\!\begin{array}{c}j\\ k
\end{array}\!\!\Big)\hbar_{-i}^{[k]}\otimes\hbar_i^{[j-k]}\otimes \XX
^j t^j )\\
&\!\!\!=\!\!\!&\mbox{$\sum\limits_{i,j=0}^{\infty}$}
\Big(\mbox{$\sum\limits_{k=0}^j$}\Big(\!\!\begin{array}{c}j
\\ k\end{array}\!\!\Big)\hbar^{[i]}\hbar_{-i}^{[k]}\otimes \XX^i\hbar_i^{[j-k]}\otimes \XX^j\Big)
\frac{(-1)^{i+j}}{i!j!}t^{i+j}\\
&\!\!\!=\!\!\!&\mbox{$\sum\limits_{i,j=0}^{\infty}$}
\Big(\mbox{$\sum\limits_{k=0}^j$}\Big(\!\!\begin{array}{c}j\\
k\end{array}\!\!\Big)\hbar^{[i+k]}\otimes\hbar^{[j-k]}\XX^i\otimes
\XX^j\Big) \frac{(-1)^{i+j}}{i!j!}t^{i+j},
\end{eqnarray*}
and
\begin{eqnarray*}
(1\otimes\DD)(Id\otimes\D^o)(\DD
)&\!\!\!=\!\!\!&\Big(\mbox{$\sum\limits_{p=0}^{\infty}$}\frac{(-1)^p
}{p!}\otimes\hbar^{[p]}\otimes
\XX^pt^p\Big)\!\cdot\!\Big(\mbox{$\sum\limits_{q=0}^{\infty}$}\frac{(-1)^q}{q!}\hbar^{[q]}
\otimes\mbox{$\sum\limits_{r=0}^q$} \Big(\!\!\begin{array}{c}q\\ r
\end{array}\!\!\Big)\XX^r\otimes \XX^{q-r}t^q\Big)\\
&\!\!\!=\!\!\!&
\mbox{$\sum\limits_{p,q=0}^{\infty}$}\Big(\mbox{$\sum\limits_{r=0}^q$}\Big(\!\!
\begin{array}{c}q\\ r\end{array}\!\!\Big)\hbar^{[q]}\otimes \hbar^{[p]}\XX^r\otimes \XX
^{p+q-r}\Big)\frac{(-1)^{p+q}}{p!q!} t^{p+q}.
\end{eqnarray*}
It is sufficient to show that the following equation holds for any fixed $m\in\Z$,
\begin{eqnarray*}
&&\mbox{$\sum\limits_{i+j=m}^{\infty}$}
\Big(\mbox{$\sum\limits_{k=0}^j$}\Big(\!\!\begin{array}{c}j\nonumber\\
k\end{array}\!\!\Big)\hbar^{[i+k]}\otimes\hbar^{[j-k]}\XX^i\otimes
\XX^j\Big)
\frac{1}{i!j!}\\
&&=\mbox{$\sum\limits_{p+q=m}^{\infty}$}\Big(\mbox{$\sum\limits_{r=0}^q$}\Big(\!\!
\begin{array}{c}q\\ r\end{array}\!\!\Big)\hbar^{[q]}\otimes \hbar^{[p]}\XX^r\otimes \XX
^{p+q-r}\Big)\frac{1}{p!q!}.
\end{eqnarray*}

Fixing $p,q,r$ such that $p+q=m,\,0\leq r\leq q$. Set $i=r,\,i+k=q$.
Then $j=m-r,\,j-k=p$. It is easy to see that the coefficients of
coefficients of $\hbar^{[q]}\otimes\hbar^{[p]}\XX^r\otimes\XX^{m-r}$
in both sides are equal to each other.\QED

\begin{lemm}\label{lemm2.5} If $\XX=L_{n_0}$, then for any $\bz \in \C,\ L_n\in\mathcal {W}$, we
have the following identities:
\begin{eqnarray}
(1 \otimes L_n) \CC_\bz\!\!\! &=&\!\!\! \mbox{$\sum\limits_{k=0}
^{\infty}$}(-1)^kb_k \CC_{\bz+k}
(\hbar_\bz^{<k>}\otimes L_{n+kn_0}t^k),\label{lemm2132} \\
(1 \otimes W_n) \CC_\bz\!\!\! &=&\!\!\! \mbox{$\sum\limits_{k=0}
^{\infty}$}(-1)^kb_k \CC_{\bz+k}
(\hbar_\bz^{<k>}\otimes W_{n+kn_0}t^k),\label{lemm201301}\\
(L_n\otimes 1 )\CC_\bz\!\!\! &=&\!\!\!\CC_{\bz-\nz} (L_n\otimes 1),\
\,L_n\mathcal{U}_\bz=\mathcal{U}_{\bz+\nz}\mbox{$\sum\limits_{k=0}
^{\infty}$}b_kL_{n+kn_0} \hbar_{-\bz+k}^{<k>}t^k,
\label{lemm2131} \\
(W_n\otimes 1 )\CC_\bz\!\!\!&=&\!\!\!\CC_{\bz-\nz} (W_n\otimes 1),\
\ W_n\mathcal{U}_\bz=\mathcal{U}_{\bz+\nz}\mbox{$\sum\limits_{k=0}
^{\infty}$}b_kW_{n+kn_0}\hbar_{-\bz+k}^{<k>}t^k,\label{lemm201302}
\end{eqnarray}
\noindent where $\nz=\frac{n}{n_0},\,b_k =
\frac{1}{k\,!}\mbox{$\prod\limits_{p=0}^{k-1}$}\big((1-p)\nz-n\big),\
b_0=1$.
\end{lemm}
\ni{\it Proof}\ \ \ Recalling (\ref{xbjxckm}), one has
\begin{eqnarray*}
&&\!\!\!\!\!\!\!\!\!(1\otimes
L_n)\CC_\bz=\mbox{$\sum\limits_{i=0}^{\infty}$}\frac{1}{i\,!}\hbar_\bz^{<i>}
\otimes \Big(\mbox{$\sum\limits_{k=0}^{i}(-1)^{k}$}\Big(
\Big(\!\!\begin{array}{c}i\\ k \end{array}\!\!\Big)
\mbox{$\prod\limits_{p=0}^{k-1}$}\big((1-p)n_0-n\big)
\XX^{i-k}L_{n+kn_0}\Big)
t^k\nonumber\\
\!\!\!&=&\!\!\!\mbox{$\sum\limits_{i=0}^{\infty}$}\,
\mbox{$\sum\limits_{k=0}^{i}$}\frac{(-1)^{k}}{k\,!(i-k)\,!}\hbar_\bz^{<i>}
\otimes \mbox{$\prod\limits_{p=0}^{k-1}$}\big((1-p)n_0-n\big)
\XX^{i-k}L_{n+kn_0}t^k\nonumber\\
\!\!\!&=&\!\!\!\mbox{$\sum\limits_{i=0}^{\infty}$}\,
\mbox{$\sum\limits_{k=0}^{\infty}$}\frac{(-1)^{k}}{k\,!i\,!}\hbar_\bz^{<i+k>}
\otimes \mbox{$\prod\limits_{p=0}^{k-1}$}\big((1-p)n_0-n\big)
\XX^{i-k}L_{n+kn_0}t^{i+k}\nonumber\\
\!\!\!&=&\!\!\!\mbox{$\sum\limits_{k=0}^{\infty}$}(-1)^{k} b_k
\mbox{$\sum\limits_{i=0}^{\infty}$}\big(\frac{1}{i\,!}\hbar_{\bz+k}^{<i>}
\otimes \XX^{i}t^i\big)\big(\hbar_\bz^{<k>} \otimes
L_{n+kn_0}t^k\big)=\mbox{$\sum\limits_{k=0}^{\infty}$}(-1)^{k} b_k
\CC_{\bz+k}\big(\hbar_\bz^{<k>} \otimes L_{n+kn_0}t^k\big).
\end{eqnarray*}
So (\ref{lemm2132}) holds. Similarly, (\ref{lemm201301}) follows by
(\ref{0801291}). From (\ref{xbjpamj}), we can deduce
\begin{eqnarray}
\!\!\!\!(L_n\!\otimes\!1)\CC_\bz\!=\!\mbox{$\sum\limits_{k=0}^{\infty}$}
\frac{1}{k\,!}L_n\hbar_\bz ^{<k>} \otimes\XX^k
t^k\!=\!\mbox{$\sum\limits_{k=0}^{\infty}$}\frac{1}{k\,!}\hbar_{\bz-\nz}
^{<k>}L_n \!\otimes\!\XX^k t^k\!=\!\CC_{\bz-\nz}(L_n\otimes1),\\
\!\!\!\!(W_n\!\otimes\!1)\CC_\bz\!=\!\mbox{$\sum\limits_{k=0}^{\infty}$}
\frac{1}{k\,!}W_n\hbar_\bz ^{<k>} \otimes\XX^k
t^k\!=\!\mbox{$\sum\limits_{k=0}^{\infty}$}\frac{1}{k\,!}\hbar_{\bz-\nz}
^{<k>}W_n \!\otimes\!\XX^k t^k\!=\!\CC_{\bz-\nz}(W_n\otimes1).
\end{eqnarray}
The first identities of (\ref{lemm2131}) and (\ref{lemm201302})
follow. As for the latter part of (\ref{lemm2131}), one has
\begin{eqnarray*}
&&\!\!\!\!\!\!\!\!\!L_n\mathcal{U}_\bz=\mbox{$\sum\limits_{p=0}^{\infty}$}
\frac{(-1)^p}{p\,!}L_n\hbar_{-\bz}^{[p]}\XX^pt^p
=\mbox{$\sum\limits_{p=0}^{\infty}$}
\frac{(-1)^p}{p\,!}\hbar_{-\nz-\bz}^{[p]}L_n\XX^p t^p\\
\!\!\!&=&\!\!\! \mbox{$\sum\limits_{p=0}^{\infty}$}
\frac{(-1)^p}{p\,!}\hbar_{-\nz-\bz}^{[p]}\Big(\mbox{$\sum\limits_{k=0}^{p}(-1)^{k}$}
\Big(\!\!\begin{array}{c}p\\k\end{array}\!\!\Big)
\mbox{$\prod\limits_{q=0}^{k-1}$}\big((1-p)n_0-n\big)
\XX^{p-k}L_{n+kn_0}\Big)t^p\\
\!\!\!&=&\!\!\! \mbox{$\sum\limits_{p=0}^{\infty}$}
\frac{(-1)^p}{(p-k)\,!}\hbar_{-\nz-\bz}^{[p]}\mbox{$\sum\limits_{k=0}^{p}(-1)^{k}$}b_k
\XX^{p-k}L_{n+kn_0}t^p\\
\!\!\!&=&\!\!\! \mbox{$\sum\limits_{k=0}^{\infty}$}\,
\mbox{$\sum\limits_{p=0}^{\infty}$}
\big(\frac{(-1)^p}{p\,!}\hbar_{-\nz-\bz}^{[p]}\XX^{p}t^{p}\big)\hbar_{-\nz-\bz}^{[k]}
b_kL_{n+kn_0}t^{k}\\
\!\!\!&=&\!\!\!\mathcal{U}_{\nz+\bz}\mbox{$\sum\limits_{k=0}^{\infty}$}\hbar_{-\nz-\bz}^{[k]}
b_kL_{n+kn_0}t^{k}=\mathcal{U}_{\nz+\bz}\mbox{$\sum\limits_{k=0}^{\infty}$}
b_kL_{n+kn_0}\hbar_{k-\bz}^{[k]}t^{k}.
\end{eqnarray*}
The second identity of (\ref{lemm201302}) can be similarly obtained.
Then the lemma follows.\QED

\ni{\it {Proof of Theorem \ref{theo1.4}}}\ \ \ For any $L_n,\,W_n\in
\mathcal {W}$, one has
\begin{eqnarray*}
\D(L_n)\!\!\! &=&\!\!\!\ \DD(L_n\otimes 1)\DD^{-1}+\DD(1\otimes
L_n)\DD^{-1}=
\DD(L_n\otimes 1)\CC+\DD(1\otimes L_n)\CC\\
\!\!\!&=&\!\!\! \DD
\CC_{-\nz}(L_n\otimes1)+\DD\mbox{$\sum\limits_{k=0}
^{\infty}$}(-1)^kb_k \CC_{k}
(\hbar^{<k>}\otimes L_{n+kn_0}t^k)\\
\!\!\!&=&\!\!\! \big(1 \otimes(1-\XX t)^\nz\big)(L_n\otimes 1)
+\mbox{$\sum\limits_{k=0} ^{\infty}$}(-1)^kb_k\big(1\otimes(1-\XX
t)^{-k}\big)
(\hbar^{<k>}\otimes L_{n+kn_0}t^k)\\
\!\!\!&=&\!\!\! L_n\otimes(1-\XX t)^\nz +\mbox{$\sum\limits_{k=0}
^{\infty}$}(-1)^kb_k\hbar^{<k>}\otimes (1-\XX
t)^{-k}L_{n+kn_0}t^k,\\
\D(W_n)\!\!\! &=&\!\!\!\ \DD(W_n\otimes 1)\DD^{-1}+\DD(1\otimes
W_n)\DD^{-1}=
\DD(W_n\otimes 1)\CC+\DD(1\otimes W_n)\CC\\
\!\!\!&=&\!\!\! \DD
\CC_{-\nz}(W_n\otimes1)+\DD\mbox{$\sum\limits_{k=0}
^{\infty}$}(-1)^kb_k \CC_{k} (\hbar^{<k>}\otimes W_{n+kn_0}t^k)\\
\!\!\!&=&\!\!\! \big(1 \otimes(1-\XX t)^\nz\big)(W_n\otimes 1)
+\mbox{$\sum\limits_{k=0} ^{\infty}$}(-1)^kb_k\big(1\otimes(1-\XX
t)^{-k}\big)
(\hbar^{<k>}\otimes W_{n+kn_0}t^k)\\
\!\!\!&=&\!\!\! W_n\otimes(1-\XX t)^\nz +\mbox{$\sum\limits_{k=0}
^{\infty}$}(-1)^kb_k\hbar^{<k>}\otimes (1-\XX t)^{-k}W_{n+kn_0}t^k,\\
\mathcal {S}(L_n)\!\!\!&=&\!\!\! \mathcal{U}^{-1}\mathcal
{S}^o(L_n)\mathcal{U}= -\mathcal{V}
\mathcal{U}_{\nz}\mbox{$\sum\limits_{k=0} ^{\infty}$}b_kL_{n+kn_0}
\hbar_{k}^{<k>}t^k=-(1-\XX t)^{-\nz}\mbox{$\sum\limits_{k=0}
^{\infty}$}b_kL_{n+kn_0} \hbar_{k}^{<k>}t^k,
\end{eqnarray*}
\begin{eqnarray*}
\mathcal {S}(W_n)\!\!\!&=&\!\!\! \mathcal{U}^{-1}\mathcal
{S}^o(W_n)\mathcal{U}= -\mathcal{V}
\mathcal{U}_{\nz}\mbox{$\sum\limits_{k=0} ^{\infty}$}b_kW_{n+kn_0}
\hbar_{k}^{<k>}t^k =-(1-\XX t)^{-\nz}\mbox{$\sum\limits_{k=0}
^{\infty}$}b_kW_{n+kn_0} \hbar_{k}^{<k>}t^k.
\end{eqnarray*}
By now, we have completed the proof of Theorem
\ref{theo1.4}.\QED\vskip10pt

\begin{lemm}\label{lemm2.6} If $\XX=W_{n_0}$, then for any $\bz \in \C,\ L_n\in\mathcal {W}$, we
have the following identities:
\begin{eqnarray}
(L_n\otimes 1 )\CC_\bz\!\!\! &=&\!\!\!\CC_{\bz-\nz} (L_n\otimes 1),\
\,(W_n\otimes 1 )\CC_\bz=\CC_{\bz-\nz} (W_n\otimes
1),\label{lemm260}\\
(1 \otimes W_n) \CC_\bz\!\!\!&=&\!\!\!\CC_\bz(1 \otimes W_n),\ \
L_n\mathcal{U}_\bz=\mathcal{U}_{\bz+\nz}\mbox{$\sum\limits_{k=0}
^{\infty}$}b_kW_{n+kn_0}
\hbar_{-\bz+k}^{<k>}t^k,\label{lemm261}\\
(1 \otimes L_n) \CC_\bz\!\!\! &=&\!\!\! \mbox{$\sum\limits_{k=0}
^{\infty}$}(-1)^kb_k \CC_{\bz+k} (\hbar_\bz^{<k>}\otimes
W_{n+kn_0}t^k),\ \
W_n\mathcal{U}_\bz=\mathcal{U}_{\bz+\nz}W_{n}.\label{lemm262}
\end{eqnarray}
\end{lemm}
\ni{\it Proof}\ \ \ Both identities of (\ref{lemm260}) can be
obtained similar to those given in Lemma \ref{lemm2.5}. The former
identity of (\ref{lemm261}) follows from
\begin{eqnarray*}
(1\otimes
W_n)\CC_\bz=\mbox{$\sum\limits_{i=0}^{\infty}$}\frac{1}{i\,!}\hbar_\bz^{<i>}
\otimes
W_{n}\XX^{i}t^k=\mbox{$\sum\limits_{i=0}^{\infty}$}\frac{1}{i\,!}\hbar_\bz^{<i>}
\otimes \XX^{i}W_{n}t^k=\CC_\bz(1\otimes W_n).
\end{eqnarray*}
As for the latter part of (\ref{lemm261}), one has
\begin{eqnarray*}
&&\!\!\!\!\!\!\!\!\!L_n\mathcal{U}_\bz=\mbox{$\sum\limits_{p=0}^{\infty}$}
\frac{(-1)^p}{p\,!}L_n\hbar_{-\bz}^{[p]}\XX^pt^p
=\mbox{$\sum\limits_{p=0}^{\infty}$}
\frac{(-1)^p}{p\,!}\hbar_{-\nz-\bz}^{[p]}L_n\XX^p t^p\\
\!\!\!&=&\!\!\! \mbox{$\sum\limits_{p=0}^{\infty}$}
\frac{(-1)^p}{p\,!}\hbar_{-\nz-\bz}^{[p]}\Big(\mbox{$\sum\limits_{k=0}^{p}(-1)^{k}$}
\Big(\!\!\begin{array}{c}p\\k\end{array}\!\!\Big)
\mbox{$\prod\limits_{q=0}^{k-1}$}\big((1-p)n_0-n\big)
\XX^{p-k}W_{n+kn_0}\Big)t^p\\
\!\!\!&=&\!\!\! \mbox{$\sum\limits_{p=0}^{\infty}$}
\frac{(-1)^p}{(p-k)\,!}\hbar_{-\nz-\bz}^{[p]}\mbox{$\sum\limits_{k=0}^{p}(-1)^{k}$}b_k
\XX^{p-k}W_{n+kn_0}t^p\\
\!\!\!&=&\!\!\! \mbox{$\sum\limits_{k=0}^{\infty}$}\,
\mbox{$\sum\limits_{p=0}^{\infty}$}
\big(\frac{(-1)^p}{p\,!}\hbar_{-\nz-\bz}^{[p]}\XX^{p}t^{p}\big)\hbar_{-\nz-\bz}^{[k]}
b_kW_{n+kn_0}t^{k}\\
\!\!\!&=&\!\!\!\mathcal{U}_{\nz+\bz}\mbox{$\sum\limits_{k=0}^{\infty}$}\hbar_{-\nz-\bz}^{[k]}
b_kW_{n+kn_0}t^{k}=\mathcal{U}_{\nz+\bz}\mbox{$\sum\limits_{k=0}^{\infty}$}
b_kW_{n+kn_0}\hbar_{k-\bz}^{[k]}t^{k}.
\end{eqnarray*}
Recalling (\ref{0801291}), one has
\begin{eqnarray*}
&&\!\!\!\!\!\!\!\!\!(1\otimes
L_n)\CC_\bz=\mbox{$\sum\limits_{i=0}^{\infty}$}\frac{1}{i\,!}\hbar_\bz^{<i>}
\otimes \Big(\mbox{$\sum\limits_{k=0}^{i}(-1)^{k}$}\Big(
\Big(\!\!\begin{array}{c}i\\ k \end{array}\!\!\Big)
\mbox{$\prod\limits_{p=0}^{k-1}$}\big((1-p)n_0-n\big)
\XX^{i-k}W_{n+kn_0}\Big)
t^k\nonumber\\
\!\!\!&=&\!\!\!\mbox{$\sum\limits_{i=0}^{\infty}$}\,
\mbox{$\sum\limits_{k=0}^{i}$}\frac{(-1)^{k}}{k\,!(i-k)\,!}\hbar_\bz^{<i>}
\otimes \mbox{$\prod\limits_{p=0}^{k-1}$}\big((1-p)n_0-n\big)
\XX^{i-k}W_{n+kn_0}t^k\nonumber\\
\!\!\!&=&\!\!\!\mbox{$\sum\limits_{i=0}^{\infty}$}\,
\mbox{$\sum\limits_{k=0}^{\infty}$}\frac{(-1)^{k}}{k\,!i\,!}\hbar_\bz^{<i+k>}
\otimes \mbox{$\prod\limits_{p=0}^{k-1}$}\big((1-p)n_0-n\big)
\XX^{i-k}W_{n+kn_0}t^{i+k}\nonumber\\
\!\!\!&=&\!\!\!\mbox{$\sum\limits_{k=0}^{\infty}$}(-1)^{k} b_k
\mbox{$\sum\limits_{i=0}^{\infty}$}\big(\frac{1}{i\,!}\hbar_{\bz+k}^{<i>}
\otimes \XX^{i}t^i\big)\big(\hbar_\bz^{<k>} \otimes
W_{n+kn_0}t^k\big)
\nonumber\\
\!\!\!&=&\!\!\!\mbox{$\sum\limits_{k=0}^{\infty}$}(-1)^{k} b_k
\CC_{\bz+k}\big(\hbar_\bz^{<k>} \otimes W_{n+kn_0}t^k\big).
\end{eqnarray*}
So the former identity of (\ref{lemm262}) holds. Observing
\begin{eqnarray*}
&&\!\!\!\!\!\!\!\!\!W_n\mathcal{U}_\bz=\mbox{$\sum\limits_{p=0}^{\infty}$}
\frac{(-1)^p}{p\,!}W_n\hbar_{-\bz}^{[p]}\XX^pt^p
=\mbox{$\sum\limits_{p=0}^{\infty}$}
\frac{(-1)^p}{p\,!}\hbar_{-\nz-\bz}^{[p]}W_n\XX^p
t^p=\mathcal{U}_{\nz+\bz}W_{n},
\end{eqnarray*}
we obtain the second identity of (\ref{lemm262}). Then the lemma
follows.\QED

\ni{\it {Proof of Theorem \ref{theo1.5}}}\ \ \ For any $L_n,\,W_n\in
\mathcal {W}$, one has
\begin{eqnarray*}
\D(L_n)\!\!\! &=&\!\!\!\ \DD(L_n\otimes 1)\DD^{-1}+\DD(1\otimes
L_n)\DD^{-1}=
\DD(L_n\otimes 1)\CC+\DD(1\otimes L_n)\CC\\
\!\!\!&=&\!\!\! \DD
\CC_{-\nz}(L_n\otimes1)+\DD\mbox{$\sum\limits_{k=0}
^{\infty}$}(-1)^kb_k \CC_{k}
(\hbar^{<k>}\otimes W_{n+kn_0}t^k)\\
\!\!\!&=&\!\!\! \big(1 \otimes(1-\XX t)^\nz\big)(L_n\otimes 1)
+\mbox{$\sum\limits_{k=0} ^{\infty}$}(-1)^kb_k\big(1\otimes(1-\XX
t)^{-k}\big)
(\hbar^{<k>}\otimes W_{n+kn_0}t^k)\\
\!\!\!&=&\!\!\! L_n\otimes(1-\XX t)^\nz +\mbox{$\sum\limits_{k=0}
^{\infty}$}(-1)^kb_k\hbar^{<k>}\otimes (1-\XX
t)^{-k}W_{n+kn_0}t^k,\\
\D(W_n)\!\!\! &=&\!\!\!\ \DD(W_n\otimes 1)\DD^{-1}+\DD(1\otimes
W_n)\DD^{-1}=
\DD(W_n\otimes 1)\CC+\DD(1\otimes W_n)\CC\\
\!\!\!&=&\!\!\! \DD
\CC_{-\nz}(W_n\otimes1)+\DD\CC(1\otimes W_n)\\
\!\!\!&=&\!\!\! \big(1 \otimes(1-\XX t)^\nz\big)(W_n\otimes 1)
+1\otimes W_n\\
\!\!\!&=&\!\!\! W_n\otimes(1-\XX t)^\nz +1\otimes W_n,\\
\mathcal {S}(L_n)\!\!\!&=&\!\!\! \mathcal{U}^{-1}\mathcal
{S}^o(L_n)\mathcal{U}= -\mathcal{V}
\mathcal{U}_{\nz}\mbox{$\sum\limits_{k=0} ^{\infty}$}b_kW_{n+kn_0}
\hbar_{k}^{<k>}t^k=-(1-\XX t)^{-\nz}\mbox{$\sum\limits_{k=0}
^{\infty}$}b_kW_{n+kn_0}
\hbar_{k}^{<k>}t^k,\\
\mathcal {S}(W_n)\!\!\!&=&\!\!\! \mathcal{U}^{-1}\mathcal
{S}^o(W_n)\mathcal{U}= -\mathcal{V} \mathcal{U}_{\nz}W_n=-(1-\XX
t)^{-\nz}W_n.
\end{eqnarray*}
By now, we have completed the proof of Theorem
\ref{theo1.4}.\QED
\end{document}